\documentclass[reqno,12pt,a4paper]{amsart}

\usepackage{mathrsfs}
\usepackage{amssymb}
\usepackage{amsfonts}
\usepackage{latexsym}
\usepackage{amsthm}

\usepackage{graphicx}
\def\lb{\label}

\newcommand{\er}[1]{\textrm{(\ref{#1})}}

\begin{document}


\renewcommand{\theequation}{\arabic{section}.\arabic{equation}}
\theoremstyle{plain}
\newtheorem{theorem}{\bf Theorem}[section]
\newtheorem{lemma}[theorem]{\bf Lemma}
\newtheorem{corollary}[theorem]{\bf Corollary}
\newtheorem{proposition}[theorem]{\bf Proposition}
\newtheorem{definition}[theorem]{\bf Definition}
\newtheorem{remark}[theorem]{\it Remark}

\def\a{\alpha}  \def\cA{{\mathcal A}}     \def\bA{{\bf A}}  \def\mA{{\mathscr A}}
\def\b{\beta}   \def\cB{{\mathcal B}}     \def\bB{{\bf B}}  \def\mB{{\mathscr B}}
\def\g{\gamma}  \def\cC{{\mathcal C}}     \def\bC{{\bf C}}  \def\mC{{\mathscr C}}
\def\G{\Gamma}  \def\cD{{\mathcal D}}     \def\bD{{\bf D}}  \def\mD{{\mathscr D}}
\def\d{\delta}  \def\cE{{\mathcal E}}     \def\bE{{\bf E}}  \def\mE{{\mathscr E}}
\def\D{\Delta}  \def\cF{{\mathcal F}}     \def\bF{{\bf F}}  \def\mF{{\mathscr F}}
\def\c{\chi}    \def\cG{{\mathcal G}}     \def\bG{{\bf G}}  \def\mG{{\mathscr G}}
\def\z{\zeta}   \def\cH{{\mathcal H}}     \def\bH{{\bf H}}  \def\mH{{\mathscr H}}
\def\e{\eta}    \def\cI{{\mathcal I}}     \def\bI{{\bf I}}  \def\mI{{\mathscr I}}
\def\p{\psi}    \def\cJ{{\mathcal J}}     \def\bJ{{\bf J}}  \def\mJ{{\mathscr J}}
\def\vT{\Theta} \def\cK{{\mathcal K}}     \def\bK{{\bf K}}  \def\mK{{\mathscr K}}
\def\k{\kappa}  \def\cL{{\mathcal L}}     \def\bL{{\bf L}}  \def\mL{{\mathscr L}}
\def\l{\lambda} \def\cM{{\mathcal M}}     \def\bM{{\bf M}}  \def\mM{{\mathscr M}}
\def\L{\Lambda} \def\cN{{\mathcal N}}     \def\bN{{\bf N}}  \def\mN{{\mathscr N}}
\def\m{\mu}     \def\cO{{\mathcal O}}     \def\bO{{\bf O}}  \def\mO{{\mathscr O}}
\def\n{\nu}     \def\cP{{\mathcal P}}     \def\bP{{\bf P}}  \def\mP{{\mathscr P}}
\def\r{\rho}    \def\cQ{{\mathcal Q}}     \def\bQ{{\bf Q}}  \def\mQ{{\mathscr Q}}
\def\s{\sigma}  \def\cR{{\mathcal R}}     \def\bR{{\bf R}}  \def\mR{{\mathscr R}}
\def\S{\Sigma}  \def\cS{{\mathcal S}}     \def\bS{{\bf S}}  \def\mS{{\mathscr S}}
\def\t{\tau}    \def\cT{{\mathcal T}}     \def\bT{{\bf T}}  \def\mT{{\mathscr T}}
\def\f{\phi}    \def\cU{{\mathcal U}}     \def\bU{{\bf U}}  \def\mU{{\mathscr U}}
\def\F{\Phi}    \def\cV{{\mathcal V}}     \def\bV{{\bf V}}  \def\mV{{\mathscr V}}
\def\P{\Psi}    \def\cW{{\mathcal W}}     \def\bW{{\bf W}}  \def\mW{{\mathscr W}}
\def\o{\omega}  \def\cX{{\mathcal X}}     \def\bX{{\bf X}}  \def\mX{{\mathscr X}}
\def\x{\xi}     \def\cY{{\mathcal Y}}     \def\bY{{\bf Y}}  \def\mY{{\mathscr Y}}
\def\X{\Xi}     \def\cZ{{\mathcal Z}}     \def\bZ{{\bf Z}}  \def\mZ{{\mathscr Z}}
\def\O{\Omega}

\newcommand{\gA}{\mathfrak{A}}
\newcommand{\gB}{\mathfrak{B}}
\newcommand{\gC}{\mathfrak{C}}
\newcommand{\gD}{\mathfrak{D}}
\newcommand{\gE}{\mathfrak{E}}
\newcommand{\gF}{\mathfrak{F}}
\newcommand{\gG}{\mathfrak{G}}
\newcommand{\gH}{\mathfrak{H}}
\newcommand{\gI}{\mathfrak{I}}
\newcommand{\gJ}{\mathfrak{J}}
\newcommand{\gK}{\mathfrak{K}}
\newcommand{\gL}{\mathfrak{L}}
\newcommand{\gM}{\mathfrak{M}}
\newcommand{\gN}{\mathfrak{N}}
\newcommand{\gO}{\mathfrak{O}}
\newcommand{\gP}{\mathfrak{P}}
\newcommand{\gQ}{\mathfrak{Q}}
\newcommand{\gR}{\mathfrak{R}}
\newcommand{\gS}{\mathfrak{S}}
\newcommand{\gT}{\mathfrak{T}}
\newcommand{\gU}{\mathfrak{U}}
\newcommand{\gV}{\mathfrak{V}}
\newcommand{\gW}{\mathfrak{W}}
\newcommand{\gX}{\mathfrak{X}}
\newcommand{\gY}{\mathfrak{Y}}
\newcommand{\gZ}{\mathfrak{Z}}

\def\ve{\varepsilon}   \def\vt{\vartheta}    \def\vp{\varphi}    \def\vk{\varkappa}

\def\Z{{\mathbb Z}}    \def\R{{\mathbb R}}   \def\C{{\mathbb C}}    \def\K{{\mathbb K}}
\def\T{{\mathbb T}}    \def\N{{\mathbb N}}   \def\dD{{\mathbb D}}


\def\la{\leftarrow}              \def\ra{\rightarrow}            \def\Ra{\Rightarrow}
\def\ua{\uparrow}                \def\da{\downarrow}
\def\lra{\leftrightarrow}        \def\Lra{\Leftrightarrow}


\def\lt{\biggl}                  \def\rt{\biggr}
\def\ol{\overline}               \def\wt{\widetilde}
\def\no{\noindent}


\let\ge\geqslant                 \let\le\leqslant
\def\lan{\langle}                \def\ran{\rangle}
\def\/{\over}                    \def\iy{\infty}
\def\sm{\setminus}               \def\es{\emptyset}
\def\ss{\subset}                 \def\ts{\times}
\def\pa{\partial}                \def\os{\oplus}
\def\om{\ominus}                 \def\ev{\equiv}
\def\iint{\int\!\!\!\int}        \def\iintt{\mathop{\int\!\!\int\!\!\dots\!\!\int}\limits}
\def\el2{\ell^{\,2}}             \def\1{1\!\!1}
\def\sh{\sharp}
\def\wh{\widehat}
\def\bs{\backslash}

\def\all{\mathop{\mathrm{all}}\nolimits}
\def\Area{\mathop{\mathrm{Area}}\nolimits}
\def\arg{\mathop{\mathrm{arg}}\nolimits}
\def\const{\mathop{\mathrm{const}}\nolimits}
\def\det{\mathop{\mathrm{det}}\nolimits}
\def\diag{\mathop{\mathrm{diag}}\nolimits}
\def\diam{\mathop{\mathrm{diam}}\nolimits}
\def\dim{\mathop{\mathrm{dim}}\nolimits}
\def\dist{\mathop{\mathrm{dist}}\nolimits}
\def\Im{\mathop{\mathrm{Im}}\nolimits}
\def\Iso{\mathop{\mathrm{Iso}}\nolimits}
\def\Ker{\mathop{\mathrm{Ker}}\nolimits}
\def\Lip{\mathop{\mathrm{Lip}}\nolimits}
\def\rank{\mathop{\mathrm{rank}}\limits}
\def\Ran{\mathop{\mathrm{Ran}}\nolimits}
\def\Re{\mathop{\mathrm{Re}}\nolimits}
\def\Res{\mathop{\mathrm{Res}}\nolimits}
\def\res{\mathop{\mathrm{res}}\limits}
\def\sign{\mathop{\mathrm{sign}}\nolimits}
\def\span{\mathop{\mathrm{span}}\nolimits}
\def\supp{\mathop{\mathrm{supp}}\nolimits}
\def\Tr{\mathop{\mathrm{Tr}}\nolimits}
\def\BBox{\hspace{1mm}\vrule height6pt width5.5pt depth0pt \hspace{6pt}}
\def\where{\mathop{\mathrm{where}}\nolimits}
\def\as{\mathop{\mathrm{as}}\nolimits}


\newcommand\nh[2]{\widehat{#1}\vphantom{#1}^{(#2)}}
\def\dia{\diamond}

\def\Oplus{\bigoplus\nolimits}



\def\qqq{\qquad}
\def\qq{\quad}
\let\ge\geqslant
\let\le\leqslant
\let\geq\geqslant
\let\leq\leqslant
\newcommand{\ca}{\begin{cases}}
\newcommand{\ac}{\end{cases}}
\newcommand{\ma}{\begin{pmatrix}}
\newcommand{\am}{\end{pmatrix}}
\renewcommand{\[}{\begin{equation}}
\renewcommand{\]}{\end{equation}}
\def\eq{\begin{equation}}
\def\qe{\end{equation}}
\def\[{\begin{equation}}
\def\bu{\bullet}

\title[{Estimates of length of spectrum}]
        {Spectral estimates for $2D$ periodic Jacobi operators}
\date{\today}

\def\Wr{\mathop{\rm Wr}\nolimits}
\def\BBox{\hspace{1mm}\vrule height6pt width5.5pt depth0pt \hspace{6pt}}

\def\Diag{\mathop{\rm Diag}\nolimits}

\date{\today}
\author[Anton Kutsenko]{Anton Kutsenko}
\address{Laboratoire de M\'ecanique Physique, UMR CNRS 5469,
Universit\'e Bordeaux 1, Talence 33405, France,  \qqq email \
kucenkoa@rambler.ru }


\maketitle

\begin{abstract}
For $2D$ periodic Jacobi operators we obtain the estimate of the
Lebesgue measure of the spectrum and estimates of edges of spectral
bands.
\end{abstract}

\section{Introduction}
\setcounter{equation}{0}

We consider a self-adjoint $2D$ periodic Jacobi operator $
\cJ:\ell^2(\Z\ts\Z)\to\ell^2(\Z\ts\Z)$ given by
\[\lb{001}
 (\cJ f)_n={A}_n f_{n+1}+B_ny_n+{A}_{n-1}^*f_{n-1},\ \ f=(f_n)\in\ell^2(\Z^2),\ \
 f_n\in\ell^2(\Z),\ \ n\in\Z.
\]
Operators $A_n,B_n:\ell^2(\Z)\to\ell^2(\Z)$ are $1D$ periodic Jacobi
matrices given by
\[\lb{002}
 A_n=SA^0_n+A^1_n+(SA^0_n)^*,\ \
 A^j_n=\diag(a^j_{nm})_{m\in\Z},\ \ j=0,1,
\]
\[\lb{003}
 B_n=SB^0_n+B^1_n+(SB^0_n)^*,\ \
 B^j_n=\diag(b^j_{nm})_{m\in\Z},\ \ j=0,1,
\]
where $S:\ell^2(\Z)\to\ell^2(\Z)$ is a shift operator
($(Sz)_m=z_{m-1}$). Complex periodic sequences $a^j,b^j$ satisfy
\[\lb{004}
 a^j_{n+p_1,m+p_2}=a^j_{nm},\ \ b^j_{n+p_1,m+p_2}=b^j_{nm},\ \
 b^1_{nm}\in\R,\ \ \forall (n,m)\in\Z^2,\ \ j=0,1,
\]
here $p_1,p_2\in\N$ are periods, we suppose that $p_1,p_2\ge3$.

For example if $A_n=I$ (identical operator), $n\in\Z$ and
$B_n=S+S^{-1}+B^1_n$, $n\in\Z$ for some diagonal $B^1_n$ then the
corresponding operator $\cJ$ is a discrete $2D$ Schr\"oedinger
operator.

By Theorem \ref{T1} the spectrum of operator $\cJ$ is
\[\lb{eq1}
 \s(\cJ)=\bigcup_{(x,y)\in[0,2\pi]^2}\{\l_n(x,y)\}_{1}^{p_1p_2},
\]
where $\l_1\le...\le\l_{p_1p_2}$ are eigenvalues of $p_1p_2\ts
p_1p_2$ matrix $J$ \er{J}. Now we give estimates of $\l_n$ and
estimates of Lebesgue measure of spectrum $\s(\cJ)$.
\begin{theorem} \lb{T2}
i) The following estimates are fulfilled
\[\lb{eq2}
 \l_n^-\le\l_n(x,y)\le\l_n^+,\ \ n\in[1,p_1p_2],\ \
 (x,y)\in[0,2\pi]^2,
\]
where $\l_n^{\pm}$ are given by \er{ld} and does not depend on
$(x,y)$.

ii) The Lebesgue measure of the spectrum satisfy
\[\lb{eq3}
 {\rm mes}(\s(\cJ))\le
 \min_{(\a,\b)\in\Z^2}r_{\a,\b},
\]
where
\[\lb{eq4}
 r_{\a,\b}=4\sum_{n=1}^{p_1}|b_{n\b}^0|+8\sum_{n=1}^{p_1}|a_{n\b}^0|-8|a^0_{\a\b}|+
 8\sum_{m=1}^{p_2}|a_{\a m}^0|+4\sum_{m=1}^{p_2}|a_{\a
 m}^1|.
\]
\end{theorem}

{\bf Remark 1.} Since $\cJ$ is a 3-diagonal matrix then we can write
simple estimate of the Lebesgue measure of spectrum
\[\lb{s1}
 {\rm mes}(\s(\cJ))\le 2\|\cJ\|\le
 2(2\max_n\|A_n\|+\max_n\|B_n\|)\le
\]
\[\lb{s2}
 \max_{n}(8\max_m|a_{nm}^0|+4\max_m|a_{nm}^1|)+\max_n(4\max_m|b_{nm}^0|+2\max_m|b_{nm}^1|).
\]
For many cases estimate \er{eq3}-\er{eq4} is better than \er{s2}
because the first depends on fewer parameters than the second
estimate, in particular \er{eq3}-\er{eq4} does not depend on $b^1$.
For example if some elements of sequence $b^1$ are large values then
\er{eq3}-\er{eq4} is better than \er{s2}.

The second reason is that we have $\min$ in \er{eq3}-\er{eq4}
instead of $\max$ in \er{s2}, but in the first estimate there is a
sum, which gives the worst result in some cases than the second
estimate.

{\no\bf 2.} For the $1D$ scalar Jacobi operator $\cA$ we have the
estimate (see e.g. \cite{DS}, \cite{Ku}, \cite{KKr})
\[\lb{002a}
 {\rm mes}(\s(\cA))\le4|a_1a_2...a_p|^{1\/p},
\]
where $a_n$ are off-diagonal elements  and $b_n$ are diagonal
elements of $\cA$.  We reach equality for the case of discrete
Shr\"odinger operator $\cA^0$ with $a_n^0=1$, $b_n^0=0$. The result
similar to \er{002a} was obtained in \cite{PR} for general non
periodic $1D$ scalar case. Recently author obtain in \cite{Ku1}
better estimate (generalized also to the vector case, when $a_n$,
$b_n$ are finite matrices)
\[\lb{003a}
 {\rm mes}(\s(\cA))\le4\min_n|a_n|.
\]
In the present paper we apply similar technique to obtain estimates
in $2D$ case. As in $1D$ case we get that the Lebesgue measure of
spectrum does not depend on diagonal elements of Jacobi matrix, but
there is a sum of some off-diagonal elements which was not in $1D$
scalar case. Nevertheless, the presence of the sum in \er{eq4}
justified (see examples below for which estimates \er{eq3}-\er{eq4}
are sharp).

{\bf 3.} Consider the important case $A_n$, $n\in\Z$ are diagonal
operators. In particular if $A_n$, $n\in\Z$ are identical operators
and $B_n$ are discrete $1D$ Schr\"odinger operators then $\cJ$
becomes $2D$ discrete Schr\"odinger operator.

If $A_n$, $n\in\Z$ are diagonal operators then all $a^1_{nm}=0$ and
\er{eq3}-\er{eq4} become
\[\lb{sch}
 {\rm mes}(\s(\cJ))\le
 4\min_{\b}\sum_{n=1}^{p_1}|b_{n\b}^0|+
 4\min_{\a}\sum_{m=1}^{p_2}|a_{\a
 m}^1|.
\]
Note that \er{s1}-\er{s2} in this case become
\[\lb{s3}
 {\rm mes}(\s(\cJ))\le4\max_{n,m}|a_{n
 m}^1|+4\max_{n,m}|b_{nm}^0|+2\max_{n,m}|b_{nm}^1|,
\]
where the dependence of $b^1$ appears again. This means that if
$b^1$ has some large components then \er{s3} is worse than \er{sch}.
We construct examples for which estimate \er{sch} is sharp. We skip
simple examples and find such examples for which the sum of $p_1$ or
$p_2$ elements in \er{sch} plays important role.

{\bf Example 1.} Let $A_n\ev0$, $n\in\Z$ and let
$$
 B_n=S+S^{-1}+4n I,\ \ n\in[1,p_1],
$$
where $S$ (see after \er{003}) is a shift operator and $I$ is
identical operator, i.e. $B_n$ is a shifted discrete $1D$
Shr\"odinger operator. By \er{001} the spectrum of $\cJ$ is
$$
 \s(\cJ)=\bigcup_{n=1}^{p_1}\s(B_n)=\bigcup_{n=1}^{p_1}[-2+4n,2+4n]=[2,2+4p_1],
$$
i.e. ${\rm mes}(\s(\cJ))=4p_1$. Estimate \er{sch} gives us the same
result ${\rm mes}(\s(\cJ))\le4p_1$.

{\no\bf 2.} Now we suppose that all $A_n\ev I$ are identical
operators and all $B_n=\diag(4m\mod4p_2)_{m\in\Z}$. In this case
$J(x,y)$ \er{J} depends on $y$ only and has a form
\[\lb{Jex2}
 J(x,y)=J(y)=\left(\begin{array}{ccccc}
                                      \hat B & \hat I & 0 & ... & e^{iy}\hat I \\
                                      \hat I & \hat B & \hat I & ... & 0 \\
                                      0 & \hat I & \hat B & ... & 0 \\
                                      ... & ... & ... & ... & ... \\
                                      e^{-iy}\hat I & 0 & 0 & ... &
                                      \hat B
  \end{array}\right),\ \ \hat B=\diag(4m-4)_1^{p_2},\ \ \hat
  I=\diag(1)_1^{p_2}.
\]
Such $J(y)$ corresponds to the matrix-valued $1D$ Jacobi operator
$\cJ_1$ with ${\rm mes}(\s(\cJ_1))=4p_1$ (this is a direct sum of
shifted Schr\"odinger operators, for more details see Example in
\cite{Ku1}). Then ${\rm mes}(\s(\cJ))=4p_2$ and \er{sch} gives us
the same result ${\rm mes}(\s(\cJ))\le4p_2$.

\section{Direct integral}
\setcounter{equation}{0}

For any $n\in\Z$ we introduce matrices $\hat A_n\ev\hat A_n(x)$ and
$\hat B_n\ev\hat B_n(x)$, $x\in[0,2\pi]$ by the following identities
\[\lb{ab}
 \hat A_n\ev\left(\begin{array}{ccccc}
                                      a^1_{n1} & \ol{a^0_{n1}} & 0 & ... & e^{ix}a^0_{np_2} \\
                                      a^0_{n1} & a^1_{n2} & \ol{a^0_{n1}} & ... & 0 \\
                                      0 & a^0_{n2} & a^1_{n3} & ... & 0 \\
                                      ... & ... & ... & ... & ... \\
                                      e^{-ix}\ol{a^0_{np_2}} & 0 & 0 & ... & a^1_{np_2}
  \end{array}\right),\ \hat B_n\ev\left(\begin{array}{ccccc}
                                      b^1_{n1} & \ol{b^0_{n1}} & 0 & ... & e^{ix}b^0_{np_2} \\
                                      b^0_{n1} & b^1_{n2} & \ol{b^0_{n1}} & ... & 0 \\
                                      0 & b^0_{n2} & b^1_{n3} & ... & 0 \\
                                      ... & ... & ... & ... & ... \\
                                      e^{-ix}\ol{b^0_{np_2}} & 0 & 0 & ... & b^1_{np_2}
                                      \end{array}\right).
\]
Also introduce the matrix $J\ev J(x,y)$, $(x,y)\in[0,2\pi]^2$ by the
following identity
\[\lb{J}
 J\ev\left(\begin{array}{ccccc}
                                      \hat B_1 & \hat A_1^* & 0 & ... & e^{iy}\hat A_{p_1} \\
                                      \hat A_1 & \hat B_2 & \hat A_2^* & ... & 0 \\
                                      0 & \hat A_2 & \hat B_3 & ... & 0 \\
                                      ... & ... & ... & ... & ... \\
                                      e^{-iy}\hat A_{p_1}^* & 0 & 0 & ... &
                                      \hat B_{p_1}
  \end{array}\right).
\]
Now we give representation of $\cJ$ \er{001} as a direct integral of
finite matrices $J$ \er{J}. Methods of the Proof of the next Theorem
are similar to methods from \cite{RS}, XIII.16, p.279.
\begin{theorem}
\lb{T1} Operator $\cJ$ is unitarily equivalent to the operator
$\mJ=\int_{[0,2\pi]^2}^{\os}J(x,y)dxdy$ acting in
$\int_{[0,2\pi]^2}^{\os}\cH dx$, where $\cH=\C^{p_1p_2}$ and matrix
$J(x,y)$ is given by \er{J}.
\end{theorem}

{\no\it Proof.} Let unitary operator $\mF_{p_2}:\ell^2(\Z)\to
L^2_{p_2}\ev\prod_{k=1}^{p_2}L^2(0,2\pi)$ be given by
\[\lb{f001}
 \mF_{p_2}((\hat f_m)_{m\in\Z})=(f_1(x),..,f_{p_2}(x)),\ \ {\rm where}\ \
 f_k(x)\ev\sum_{m=-\iy}^{+\iy}\hat f_{k+p_2m}e^{inx},\ \ k\in[1,p_2].
\]
Let unitary operator ${\bf F}:\ell^2(\Z\ts\Z)\to{\bf
L}^2_{p_2}\ev\prod_{k\in\Z}L^2_{p_2}$ be given by
\[\lb{f002}
 {\bf F}((f_n)_{n\in\Z})=(\mF_{p_2}f_n)_{n\in\Z},\ \
 f_n\in\ell^2(\Z),\ \ n\in\Z.
\]
Introduce operator
\[\lb{fj1}
 \cJ_1\ev{\bf F}\cJ{\bf F}^{-1}:{\bf
L}^2_{p_2}\to{\bf L}^2_{p_2}.
\]
Substituting \er{f001}-\er{fj1} into \er{001}-\er{003} and using
\er{004} we deduce that $\cJ_1$ has a form
\[\lb{f003}
 (\cJ_1 f)_n=\hat{A}_n f_{n+1}+\hat B_nf_n+\hat{A}_{n-1}^*f_{n-1},\
 \ f=(f_n)\in{\bf L}^2_{p_2},\ \ f_n\in L^2_{p_2},\ \ n\in\Z
\]
where
\[\lb{f004}
 \hat{A}_n=\mF_{p_2}{A}_n\mF_{p_2}^{-1},\ \
 \hat{B}_n=\mF_{p_2}{B}_n\mF_{p_2}^{-1}
\]
is given by \er{ab}. Let unitary operator ${\bf F}_1:{\bf
L}^2_{p_2}\to{\bf L}^2\ev\prod_{k=1}^{p_1p_2} L^2([0,2\pi]^2)$ be
given by
\[\lb{f006}
 {\bf F}_1((f_n)_{n\in\Z})=(g_1,...,g_{p_1}),\ \ {\rm where}\ \
 g_k\ev\sum_{n=-\iy}^{+\iy}f_{k+p_1n}e^{iny},\ \ k\in[1,p_1],
\]
here $f_n\ev f_n(x)\in L^2_{p_2}$ for any $n\in\Z$. Introduce
operator
\[\lb{fj2}
 \cJ_2\ev{\bf F}_1\cJ_1{\bf F}_1^{-1}:{\bf L}^2\to{\bf
 L}^2.
\]
Substituting \er{f006}-\er{fj2} into \er{f003} we deduce that
$\cJ_2$ has a form
$\cJ_2((g_n)_1^{p_1p_2})=J_2(x,y)(g_n)_1^{p_1p_2}$, where $g_n\in
L^2([0,2\pi]^2)$ and matrix $J(x,y)$ is given by \er{J}, i.e.
$\cJ_2=\mJ$. By \er{fj1}, \er{fj2} we deduce that $\cJ$ unitarily
equivalent to the operator $\mJ$. \BBox

\section{Proof of Theorem \ref{T2}}
\setcounter{equation}{0}

For any $n$ let us denote
\[\lb{a0}
 \hat A_n^0\ev\left(\begin{array}{ccccc}
                                      a^1_{n1} & \ol{a^0_{n1}} & 0 & ... & 0 \\
                                      a^0_{n1} & a^1_{n2} & \ol{a^0_{n1}} & ... & 0 \\
                                      0 & a^0_{n2} & a^1_{n3} & ... & 0 \\
                                      ... & ... & ... & ... & ... \\
                                      0 & 0 & 0 & ... & a^1_{np_2}
  \end{array}\right),\ \hat A_n^1\ev\left(\begin{array}{ccccc}
                                      0 & 0 & 0 & ... & e^{ix}a^0_{np_2} \\
                                      0 & 0 & 0 & ... & 0 \\
                                      0 & 0 & 0 & ... & 0 \\
                                      ... & ... & ... & ... & ... \\
                                      e^{-ix}\ol{a^0_{np_2}} & 0 & 0 & ... &
                                      0
                                      \end{array}\right),
\]
\[\lb{b0}
 \hat B_n^0\ev\left(\begin{array}{ccccc}
                                      b^1_{n1} & \ol{b^0_{n1}} & 0 & ... & 0 \\
                                      b^0_{n1} & b^1_{n2} & \ol{b^0_{n1}} & ... & 0 \\
                                      0 & b^0_{n2} & b^1_{n3} & ... & 0 \\
                                      ... & ... & ... & ... & ... \\
                                      0 & 0 & 0 & ... & b^1_{np_2}
  \end{array}\right),\ \hat B_n^1\ev\left(\begin{array}{ccccc}
                                      0 & 0 & 0 & ... & e^{ix}b^0_{np_2} \\
                                      0 & 0 & 0 & ... & 0 \\
                                      0 & 0 & 0 & ... & 0 \\
                                      ... & ... & ... & ... & ... \\
                                      e^{-ix}\ol{b^0_{np_2}} & 0 & 0 & ... &
                                      0
                                      \end{array}\right),
\]
\[\lb{J0}
 J_0\ev\left(\begin{array}{ccccc}
                                      \hat B_1^0 & ( \hat A_1^0)^* & 0 & ... & 0 \\
                                      \hat A_1^0 & \hat B_2 &  (\hat A_2^0)^* & ... & 0 \\
                                      0 & \hat A_2^0 & \hat B_3^0 & ... & 0 \\
                                      ... & ... & ... & ... & ... \\
                                      0 & 0 & 0 & ... &
                                      \hat B_{p_1}^0
  \end{array}\right),\ J_1\ev\left(\begin{array}{ccccc}
                                      \hat B^1_1 &  (\hat A_1^1)^* & 0 & ... & e^{iy}\hat A_{p_1} \\
                                      \hat A_1^1 & \hat B_2^1 &  (\hat A_2^1)^* & ... & 0 \\
                                      0 & \hat A^1_2 & \hat B^1_3 & ... & 0 \\
                                      ... & ... & ... & ... & ... \\
                                      e^{-iy}\hat A_{p_1}^* & 0 & 0 & ... &
                                      \hat B^1_{p_1}
  \end{array}\right).
\]
Then
\[\lb{p001}
 \hat A_n=\hat A_n^0+\hat A_n^1,\ \ \hat B_n=\hat B_n^0+\hat B_n^1,\
 \ J=J_0+J_1.
\]
Note that $J_0$ does not depend on $x,y$. For any self-adjoin matrix
$A$ let us denote $A_+\ev P_+A\ge0$, $A_-\ev -P_-A\ge0$, where
$P_{\pm}$ are projectors to positive and negative eigenspaces of
$A$. Then $A=A_+-A_-$ and $|A|\ev(AA)^{\frac12}=A_++A_-$. We will
write $A\ge B$ if $A-B\ge0$, i.e. $A-B$ is positive definite
self-adjoint matrix. Now we will estimate $J_1$. Firstly define the
matrix $C$ by
\[\lb{C}
 C=\diag(C_1,...,C_{p_1}),
\]
where matrices $C_n$ is given by: for $n\in[2,p_1-1]$
\[\lb{Ck}
 C_n=|\hat B^1_{n}|+|\hat A^1_{n-1}|+|\hat
 A^1_n|=\diag(c_n,0,..,0,c_n),\ \
 c_n=|b^0_{np_2}|+|a^0_{np_2}|+|a^0_{n-1,p_2}|,
\]
for $n=1,p_1$
\[\lb{C1}
 C_1=|\hat B_1^1|+|\hat A_1^1|+D,\ \ C_{p_1}=|\hat B_{p_1}^1|+|\hat
 A_{p_1-1}^1|+D,
\]
where
\[\lb{D}
 D=\diag(|a_{p_1m}^1|)_{m=1}^{p_2}+\diag(|
 a_{p_1m}^0|)_{m=1}^{p_2}+\diag(|a_{p_1p_2}^0|,|a_{p_11}^0|,..,|a_{p_1,p_2-1}^0|).
\]
Note that $C$ is a positive definite diagonal matrix does not depend
on $x$ and $y$.
\begin{lemma} The following inequalities are fulfilled
\[\lb{l1}
 -C\le J_1\le C.
\]
\end{lemma}
{\it Proof.} Note that
\[\lb{l001}
 J_1=\diag(\hat B^1_{n})_{1}^{p_1}+\sum_{n=1}^{p_1}E_{n}+F_1+F_2+F_3,
\]
where
\[\lb{En}
 E_1=\left(\begin{array}{ccccc}
                                      0 & ( \hat A_1^0)^* & 0 & ... & 0 \\
                                      \hat A_1^0 & 0 &  0 & ... & 0 \\
                                      0 & 0 & 0 & ... & 0 \\
                                      ... & ... & ... & ... & ... \\
                                      0 & 0 & 0 & ... &  0
  \end{array}\right),\ \ E_2=\left(\begin{array}{ccccc}
                                      0 & 0 & 0 & ... & 0 \\
                                      0 & 0 &  (\hat A_2^0)^* & ... & 0 \\
                                      0 & \hat A_2^0 & 0 & ... & 0 \\
                                      ... & ... & ... & ... & ... \\
                                      0 & 0 & 0 & ... & 0
  \end{array}\right),\ ...
\]
and
\[\lb{F1}
 F_1=\left(\begin{array}{ccccc}
                                      0 & 0 & 0 & ... & e^{iy}G_1 \\
                                      0 & 0 & 0 & ... & 0 \\
                                      0 & 0 & 0 & ... & 0 \\
                                      ... & ... & ... & ... & ... \\
                                      e^{-iy}G_1^* & 0 & 0 & ... &
                                      0
                                      \end{array}\right),\ \ G_1=\left(\begin{array}{ccccc}
                                      0 & 0 & 0 & ... & e^{ix}a^0_{np_2} \\
                                      a^0_{p_11} & 0 & 0 & ... & 0 \\
                                      0 & a^0_{p_12} & 0 & ... & 0 \\
                                      ... & ... & ... & ... & ... \\
                                      0 & 0 & 0 & ... & 0
  \end{array}\right),
\]
\[\lb{F2}
 F_2=\left(\begin{array}{ccccc}
                                      0 & 0 & 0 & ... & e^{iy}G_2 \\
                                      0 & 0 & 0 & ... & 0 \\
                                      0 & 0 & 0 & ... & 0 \\
                                      ... & ... & ... & ... & ... \\
                                      e^{-iy}G_2^* & 0 & 0 & ... &
                                      0
                                      \end{array}\right),\ \ G_2=\left(\begin{array}{ccccc}
                                      0 & \ol{a^0_{p_11}} & 0 & ... & 0 \\
                                      0 & 0 & \ol{a^0_{p_11}} & ... & 0 \\
                                      0 & 0 & 0 & ... & 0 \\
                                      ... & ... & ... & ... & ... \\
                                      e^{-ix}\ol{a^0_{p_1p_2}} & 0 & 0 & ... &
                                      0
  \end{array}\right),
\]
\[\lb{F3}
  F_3=\left(\begin{array}{ccccc}
                                      0 & 0 & 0 & ... & e^{iy}G_3 \\
                                      0 & 0 & 0 & ... & 0 \\
                                      0 & 0 & 0 & ... & 0 \\
                                      ... & ... & ... & ... & ... \\
                                      e^{-iy}G_3^* & 0 & 0 & ... &
                                      0
                                      \end{array}\right),\ \
                                      G_3=\diag(a^1_{p_1n})_1^{p_2}.
\]
Using \er{C}-\er{D} and \er{En}-\er{F3} we deduce that
\[\lb{l002}
 \diag(
 |\hat B^1_{n}|)_{1}^{p_1}+\sum_{n=1}^{p_1}|E_{n}|+|F_1|+|F_2|+|F_3|=C,
\]
which with \er{l001} gives us \er{l1}. \BBox

Let us denote eigenvalues of matrices $J_0-C$ and $J_0+C$
\er{J0},\er{C} as
\[\lb{ld}
 \l_1^{-}\le...\le\l_{p_1p_2}^-\ \ {\rm and}\ \ \l_1^{+}\le...\le\l_{p_1p_2}^+.
\]

{\no\bf Proof of Theorem \ref{T2}.} i) follows from $\er{p001}_3$,
and \er{l1}, since $J_0-C\le J\le J_0+C$ and $J_0\pm C$ does not
depend on $(x,y)$.

ii) Without loss of generality, we may assume that
\[\lb{t001}
 \min_{(\a,\b)\in\Z^2}r_{\a,\b}=r_{p_1p_2}=4\sum_{n=1}^{p_1}|b_{np_2}^0|+8\sum_{n=1}^{p_1-1}|a_{np_2}^0|+
 8\sum_{m=1}^{p_2}|a_{p_1 m}^0|+4\sum_{m=1}^{p_2}|a_{p_1
 m}^1|
\]
otherwise we enumerate the periodic sequences $a^j,b^j$, $j=0,1$.
Since $\l_n^{\pm}$ does not depend on $(x,y)$ then from \er{eq2} we
obtain
\[\lb{pr1}
 {\rm mes}(\s(\cJ))\le\sum_{n=1}^{p_1p_2}(\l_n^+-\l_n^-)
\]
Now inequality \er{eq3} follows from the identity
\[\lb{t002}
 \sum_{n=1}^{p_1p_2}(\l_n^+-\l_n^-)=2\Tr C=r_{p_1p_2},
\]
where we use definition before \er{ld} and \er{C}-\er{D} with
\er{t001}. \BBox


\begin{thebibliography} {9999}
\setlength{\itemsep}{-\parskip} \footnotesize

\bibitem [CG]{CG}  Clark, S.; Gesztesy, F. On Weyl–-Titchmarsh theory for
singular finite difference Hamiltonian systems, J. Comput. Appl.
Math., 171 (2004) 151–184.

\bibitem [CGR]{CGR} Clark, S.; Gesztesy, F.; Renger, W.
Trace formulas and Borg-type theorems for matrix-valued Jacobi and
Dirac finite difference operators. J. Diff. Eq. 219 (2005),
144--182.

\bibitem[DS] {DS} P. Deift, B. Simon. Almost periodic Schr\"odinger operators III.
The absolutely continuous spectrum in one dimension. Commun. Math.
Phys., 90, 389–411 (1983).

\bibitem[K] {K} Kato T. Perturbation Theory for Linear Operators. Springer (February 15,
1995).

\bibitem [Ku] {Ku} Kutsenko A. Estimates of Parameters for Conformal Mappings Related to a Periodic Jacobi
Matrix. Journal of Mathematical Sciences,     Volume 134, Number 4 /
April 2006, Pages 2295-2304.

\bibitem [Ku1] {Ku1} Kutsenko A. Sharp spectral estimates for periodic matrix-valued Jacobi
operators, preprint 2010. http://arxiv.org/abs/1007.5412

\bibitem [KKr]{KKr}  Korotyaev, E.; Krasovsky, I.
Spectral estimates for periodic Jacobi matrices, Commun. Math. Phys.
234(2003), 517-532.


\bibitem [KKu]{KKu} Korotyaev, E., Kutsenko, A.
    Lyapunov functions for periodic matrix-valued Jacobi operators,
    AMS translations Series 2, 225 (2008), 117—-131.

\bibitem[KKu1] {KKu1} Korotyaev, E., Kutsenko, A. Borg type uniqueness Theorems for
periodic Jacobi operators with matrix valued coefficients. Proc. of
the AMS, Volume 137, Number 6, June 2009, Pages 1989–-1996.

\bibitem[L] {L} Y. Last. On the measure of gaps and spectra for discrete 1D
Schr\"odinger operators. Commun. Math. Phys., 149, 347–-360 (1992).

\bibitem[PR]{PR} A. Poltoratski, C. Remling. Reflectionless Herglotz Functions and Jacobi
Matrices, Commun. Math. Phys. Volume 288 Number 3(2009), 1007--1021.

\bibitem[RS]{RS} M. Reed ; B. Simon. Methods of modern mathematical physics. IV.
Analysis of operators. Academic Press, New York-London, 1978.

\bibitem[S]{S} B. Simon, Orthogonal polynomials on the unit circle, Part 1 and Part 2, AMS, Providence, RI, 2005.

\bibitem[Te] {Te}  G. Teschl, Jacobi Operators and Completely Integrable Nonlinear
Lattices, Mathematical Surveys and Monographs, vol. 72, AMS, Rhode
Island, 2000.

\bibitem [vM]{vM}  P. van Moerbeke. The spectrum of Jacobi matrices.
Invent. Math. 37 (1976), no. 1, 45--81.

\end{thebibliography}
\end{document}